\documentclass[12pt]{article}
\usepackage{amssymb,amsfonts,amsmath,amsthm}
\newcommand{\ol}{\overline}
\newcommand{\pp}{{\mathbb P}}
\newcommand{\Ra}{\Rightarrow}
\newcommand{\fracd}[2]{\frac {\displaystyle #1}{\displaystyle #2 }}
\newcommand{\suml}{\sum\limits}
\newcommand{\prodl}{\prod\limits}
\newcommand{\zz}{{\mathbb Z}}
\newcommand{\nn}{{\mathbb N}}
\newcommand{\ee}{{\mathbb E}}
\newcommand{\rr}{{\mathbb R}}

\newcommand{\calf}{{\mathcal F}}
\newcommand{\veps}{\varepsilon}
\newcommand{\beq}{\begin{eqnarray*}}
\newcommand{\feq}{\end{eqnarray*}}
\newcommand{\beqn}{\begin{eqnarray}}
\newcommand{\feqn}{\end{eqnarray}}
\newcommand{\as}{\mbox{a.s.}}

\newcommand{\proda}{{\mbox {\Large $\Pi$}}}
\newtheorem{theorem}{Theorem}
\makeatletter \@addtoreset{theorem}{section}\makeatother

\newtheorem{lemma}[theorem]{Lemma}
\newtheorem{assume}[theorem]{Assumption}
\newtheorem{proposition}[theorem]{Proposition}
\title{A $\log$-scale limit theorem for one-dimensional random
walks in random environments}
\author{
Alexander Roitershtein\footnote{Department of Mathematics, Technion -
IIT, Haifa 32000, Israel (e-mail: roiterst@tx.technion.ac.il).}}
\begin{document}
\maketitle
\begin{abstract}
We consider a transient one-dimensional random walk $X_n$ in
random environment having zero asymptotic speed. For a class of
non-i.i.d. environments we show that $\log X_n/\log n$ converges
in probability to a positive constant.
\end{abstract}

\section{Introduction}
\label{scale} In this note we consider a one-dimensional random
walk $X_n$ in ergodic environments and prove a log-scale result
which is consistent with limit theorems known for i.i.d. and
Markovian environments, namely that $\log X_n/\log n$ converges in
probability to a positive constant.

The environment is a stationary and ergodic sequence indexed by
the sites of $\zz$ which represents the probabilities that a
nearest-neighbor random walk $X_n$ moves to the right. The
(annealed) law of $X_n$ is obtained by averaging its quenched law
(i.e. given a fixed environment) over the set of environments.

More formally, let $\Omega=(0,1)^\zz,$ $\calf$ be its Borel
$\sigma-$algebra, and $P$ be a stationary and ergodic probability
measure on $(\Omega,\calf).$ The set $\Omega$ serves as the set of
{\em environments} for a random walk $X_n$ on $\zz.$ The {\em
random walk in the environment} $\omega=\{\omega_i\}_{i \in \zz}
\in \Omega$ is the time-homogeneous Markov chain $X=$ $\{X_n\}_{n
\in {\nn}}$ taking values on $\zz$ and governed by the {\em
quenched} law \beq
P_\omega(X_0=0)=1~~~\mbox{and}~~~P_{\omega}(X_{n+1}=j|X_n=i)=\left\{
\begin{array}{ll}
\omega_i&\mbox{if}~j=i+1,\\
1-\omega_i&\mbox{if}~j=i-1.
\end{array}
\right. \feq Let $\mathcal G$ be the cylinder $\sigma-$algebra on
$\zz^\nn,$ the path space of the walk. {\em The random walk in
random environment} (RWRE) is the process $(X,\omega)$ on the
measurable space $\left(\Omega \times \zz^\nn, {\mathcal F} \times
{\mathcal G}\right)$ with the {\em annealed} law $\pp=P \otimes
P_\omega$ defined by \beq \pp(F \times G)=\int_F
P_\omega(G)P(d\omega)=E_P\left(P_\omega(G);F\right),~~~ F \in
{\mathcal F},~ G \in {\mathcal G}. \feq  We refer the reader to
\cite[Section 2]{notes} for a detailed exposition of the
one-dimensional model.

Let $\rho_n=(1-\omega_n)/\omega_n,$ and \beqn \label{sf}
R(\omega)=1+\sum_{n=0}^{+\infty}\rho_0\rho_{-1} \cdots \rho_{-n}.
\feqn If $E_P(\log \rho_0)<0$ then (Solomon \cite{solomon} for
i.i.d. environments and Alili \cite{alili} in the general case)
$\lim_{n \to \infty} \pp(X_n=+\infty)=1$ and \beqn \label{speed}
\mbox{v}_P:=\lim_{n \to +\infty}\fracd{X_n}{n}=
\fracd{1}{2E_P(R)-1},~\pp-\as \feqn The role played by the process
$V_n=\sum_{i=0}^n \log \rho_{-i}=\log\bigl(\rho_0\rho_{-1} \cdots
\rho_{-n}\bigr)$ in the theory of one-dimensional RWRE stems from
the explicit form of harmonic functions (cf. \cite[Section
2.1]{notes}), which allows one to relate hitting times of the RWRE
to those associated with the random walk $V_n.$

The assumptions we impose below will imply that \beqn
\label{case}E_P(\log \rho_0)<0~~~\mbox{and}~~~E_P( R)=\infty,
\feqn that is the random walk is transient to the right but
$\mbox{v}_P=0.$

The limit theorem of Kesten, Kozlov, and Spitzer \cite{kks} states
that if $P$ is a product measure and $E_P(\rho^\kappa)=1$ for some
positive $\kappa,$ then, under some additional technical
conditions, an appropriately normalized sequence $X_n$ converges
to a non-degenerate (scaled stable) limit law. In particular, if
$\kappa \in (0,1),$ then $n^{-\kappa} X_n$ converges in law and
thus $\log X_n/\log n$ converges to $\kappa$ in probability.

The limit laws established in \cite{kks} were extended to
Markov-dependent environments in \cite{mars}, where the condition
$E_P(\rho_0^\kappa)=1$ is replaced by: \beqn \label{kappa-new}
\lim_{n \to \infty} \fracd{1}{n} \log E_P \left(\proda_{i=0}^{n-1}
\rho_i^\kappa\right)=0~~~\mbox{for some}~\kappa>0. \feqn

The key element of the proofs of the limit laws in \cite{kks} and
\cite{mars} is that the distribution tail of the random variable
$R$ is regularly varying. In fact, under assumptions imposed on
the environment in \cite{kks} and \cite{mars}, the limit $\lim_{t
\to \infty}t^\kappa P(R>t)$ exists and is strictly positive for
$\kappa$ given by \eqref{kappa-new}. The existence of the limit is
closely related to a renewal theory for the random walk
$V_n=\sum_{i=0}^{n-1} \log \rho_{-i}$ (cf. \cite{kesten-randeq},
see also \cite{goldie} for an alternative proof and \cite{saporta}
and \cite[Section 2.3]{mars} for Markovian extensions), and thus
the lack of an adequate renewal theory for dependent random
variables puts limitation on the scope of this approach to limit
laws for RWRE.

In this paper, assuming that the environment is strongly mixing,
the sequence $V_n/n$ satisfies the Large Deviation Principle
(LDP), and \eqref{kappa-new} holds for some $\kappa \in (0,1),$ we
prove that $\log X_n$ is asymptotic to $\kappa \log n$ by using a
different method, based on a study of the generating function of
an associated branching process and a large deviation analysis of
the tail.

We assume that the following holds for the environment $\omega:$
\begin{assume}
\label{v239} $\omega=\{\omega_n\}_{n \in \zz}$ is a stationary
sequence such that: \item[(A1)] The series $\sum_{n=1}^\infty
n^\xi \cdot \eta\left( n \right)$ converges for any $\xi >0,$
where the strong mixing coefficients $\eta_n$ are defined by \beq
\eta(n):=\sup\{|P(A \cap B)-P(A)P(B)|:A \in \sigma(\omega_i:i \leq
0),~B \in \sigma(\omega_i: i \geq n)\}. \feq   \item[(A2)] The
process $\ol{V}_n=n^{-1}\sum_{i=0}^{n-1} \log \rho_i$ satisfies
the LDP with a good rate function $J: \rr \to$ $[0,\infty],$ that
is the level sets $J^{-1}[0,a]$ are compact for any $a>0,$ and for
any Borel set $A,$ \beq -J(A^o) \leq \liminf_{n \to \infty}
\fracd{1}{n} \log \pp (\ol{V}_n \in A) \leq \limsup_{n \to \infty}
\fracd{1}{n} \log \pp (\ol{V}_n \in A) \leq -J(A^c), \feq where
$A^o$ denotes the interior of $A,$ $A^c$ the closure of $A,$ and
$J(A):=\inf_{x \in A} J(x).$ \item[(A3)] $\limsup_{n \to \infty}
\fracd{1}{n}\log E_P \left(\prod_{i=0}^{n-1} \rho_i^\lambda
\right)<\infty$ for all $\lambda \in {\mathbb R}.$ \item [(A4)]
$\limsup_{n \to \infty} \fracd{1}{n} \log E_P \left(
\prod_{i=0}^{n-1} \rho_i\right) \geq 0$ and $\limsup_{n \to
\infty} \fracd{1}{n} \log E_P \left(\prod_{i=0}^{n-1}
\rho_i^{\lambda_1} \right)<0$ for some $\lambda_1 >0.$
\end{assume}
Assumption \ref{v239} is a modification of a condition which was
introduced in the context of RWRE in \cite[Sect. 2.4]{notes}.
Assumptions {\em (A2)} and {\em (A3)} ensure that for any $\lambda
\in {\mathbb R},$  the limit $\Lambda(\lambda):=\lim_{n \to
\infty}\fracd{1}{n}\log E_P \left(\prod_{i=0}^{n-1} \rho_i^\lambda
\right)$ exists and satisfies $\Lambda(\lambda)=\sup_x\{\lambda
x-J(x)\}$ \cite[p. 135]{amir-ofer-book}. It follows from condition
{\em (A4)} that $\Lambda(1) \geq 0$ and $\Lambda(\lambda_1)<0.$ By
Jensen's inequality, $E_P(\log \rho_0) \leq
\Lambda(\lambda_1)/\lambda_1<0,$ and hence \eqref{case} holds.
Moreover, since $\Lambda (\lambda)$ is a convex function with
$\Lambda(0)=0,\Lambda(1) \geq 0,$ and taking both positive and
negative values, there exists $\kappa \in (0,1]$ such that for
each $\lambda>0,$ \beqn \label{sign} \Lambda(\lambda)~\mbox{has
the same sign as}~\lambda-\kappa. \feqn Since $J(x) \geq
\sup_\lambda\{\lambda x -\Lambda(\lambda)\}$ (cf. \cite[Theorem
4.5.10]{amir-ofer-book}), it follows that $J(0)>0$ and hence
(similarly to \cite[Sect 2.4]{notes}): \beqn \label{kappa}
\kappa=\min_{y > 0} J(y)/y. \feqn For conditions on Markov and
general stationary processes implying Assumption \ref{v239} we
refer to \cite{bryc-dembo, amir-ofer-book,guliver} and references
therein. In particular, all bounded, stationary and strongly
mixing sequences $(\rho_n)_{n \in \zz}$ with a fast mixing rate
(namely $\frac{\log \eta(n)}{n (\log n)^{1+\delta}} \to -\infty$
for some $\delta>0$) do satisfy Assumptions {\em (A1)--(A3)} (cf.
\cite{bryc-dembo}).

Our main result is:
\begin{theorem}
\label{scale-main} Suppose that Assumption \ref{v239} holds. Then,
\beq n^{-\alpha} X_n \overset{\pp}{\Ra}
\begin{cases}
\infty&\mbox{\em if}~\alpha<\kappa\\
0&\mbox{\em if}~\alpha>\kappa,
\end{cases}
\feq where $\kappa$ is defined in \eqref{kappa}.
\end{theorem}
This theorem is deduced from a similar result for the hitting time
$T_n,$ defined by \beqn \label{ti} T_n=\min\{i: X_i=n\},~~~n \in
\nn. \feqn We have:
\begin{proposition}
\label{scale-ti} Suppose that Assumption \ref{v239} holds. Then,
\beq n^{-1/\alpha} T_n \overset{\pp}{\Ra}
\begin{cases}
0&\mbox{\em if}~\alpha<\kappa\\
\infty&\mbox{\em if}~\alpha>\kappa,
\end{cases}
\feq where $\kappa$ is defined in \eqref{kappa}.
\end{proposition}
The proof of Proposition \ref{scale-ti}, included in
Section~\ref{proofs}, is based on the study of generating
functions of the branching process $Z_n.$ Once this proposition is
obtained, Theorem \ref{scale-main} is derived from it as follows
(cf. \cite{kks}). For any positive integers $\eta,\zeta,n,$ \beq
\{T_\zeta \geq n\} \subset \{ X_n \leq \zeta\} \subset
\{T_{\zeta+\eta} \geq n\} \bigcup \{\inf_{i \geq T_{\zeta+\eta} }
X_i -(\zeta+\eta) \leq -\eta\}.\feq Since the random variables
$\inf_{i \geq T_{\zeta+\eta} } X_i -(\zeta+\eta)$ and $\inf_{i
\geq 0} X_k$ have the same annealed distribution, \beqn
\label{unif-239} \pp(T_\zeta \geq n) \leq \pp( X_n \leq \zeta)
\leq  P(T_{\zeta+\eta} \geq n) + \pp(\mbox{inf}_{i \geq 0 } X_i
\leq -\eta). \feqn  Fix any numbers $\eta \in \nn,$ $x>0,\alpha
\in (0,\kappa) \cup (\kappa,\infty),$ and let $\zeta(n)=[x
n^\alpha],$ where $[t]$ denotes the integer part of $t.$ It
follows from \eqref{unif-239} that \beq && \limsup_{n \to \infty}
P( n^{-\alpha} X_n \leq x)\leq \limsup_{n \to \infty} P(X_n \leq
\zeta(n)+1) \leq \liminf_{n \to \infty} P(T_{\zeta(n)+1+\eta} \geq
n) \\&& \qquad  + \pp(\mbox{inf}_{i \geq 0 } X_i \leq -\eta)
=\lim_{\zeta \to \infty} P(T_\zeta\geq \zeta^{1/\alpha}
x^{-1/\alpha})+ \pp(\mbox{inf}_{i \geq 0 } X_i \leq -\eta). \feq
Similarly, \beq && \liminf_{n \to \infty} P( n^{-\alpha} X_n \leq
x)\geq \liminf_{n \to \infty} P(X_n \leq \zeta(n)) \geq \limsup_{n
\to \infty} P(T_{\zeta(n)} \geq n) \\&& \qquad =\lim_{\zeta \to
\infty} P(T_\zeta\geq \zeta^{1/\alpha} x^{-1/\alpha}). \feq Since
$X_n$ is transient to the right, $\pp(\mbox{inf}_{i \geq 0 } X_i
\leq -\eta)$ can be made arbitrary small by fixing $\eta$ large.
Thus, \beq && \lim_{n \to \infty} P( n^{-\alpha} X_n \leq
x)=\lim_{\zeta \to \infty} P(T_\zeta\geq \zeta^{1/\alpha}
x^{-1/\alpha})=\begin{cases}
0&\mbox{if}~\alpha<\kappa\\
1&\mbox{if}~\alpha>\kappa.
\end{cases} \feq The next
section is devoted to the proof of Proposition \ref{scale-ti}.
\section{Proof of Proposition \ref{scale-ti}}
\label{proofs} First, we consider a branching process $Z_n$ in
random environment, closely related to the RWRE. The hitting times
$T_n$ of the random walk are associated by \eqref{ti-zed} and
\eqref{ti-minus} below to the partial sums of the branching
process and Proposition \ref{scale-ti} is then derived from a
corresponding scale result for $\sum_{i=1}^n Z_i$ (see Proposition
\ref{scale-zed}).

The branching process was introduced in the context of the RWRE in
\cite{kozlov73}, and has been exploited in several works including
\cite{kks} (a detailed description of its construction can be
found e.g. in \cite{many-visits}). Let \beq U_i^n=\#
\{t<T_n:X_t=i,~X_{t+1}=i-1\},~~~n \in \nn,~i \in \zz, \feq the
number of moves to the left at site $i$ up to time $T_n.$ Then
\beqn \label{ti-zed} T_n=n+2\sum\limits_{i=-\infty}^n U_i^n. \feqn
When $U^n_n=0,U^n_{n-1}, \ldots, U^n_{n-i+1}$ and $\omega_n,
\omega_{n-1} \ldots,\omega_{n-i}$ are given, $U^n_{n-i}$ is the
sum of $U^n_{n-i+1}+1$ i.i.d. geometric random variables that take
the value $k$ with probability $\omega_{n-i}(1-\omega_{n-i})^k,$
$k=0,1,\ldots$ Since $X_n$ is transient to the right we have:
\beqn \label{ti-minus} \sum\limits_{i \leq 0} U_i^n
<\infty,~~~\pp-\as\feqn Therefore, in order to prove Proposition
\ref{scale-ti} it is sufficient to show that the corresponding
result holds for the sums $\sum_{i=1}^n U_i^n.$

These sums have the same distribution as \beq \label{zed}
\sum\limits_{i=0}^{n-1}Z_i, \feq where $Z_0=0,Z_1,Z_2,\ldots$
forms a branching process in random environment with one immigrant
at each unit of time. Without loss of generality, we shall assume
that the underlying probability space is extended to fit not only
the random walk but also the branching process. Thus, when
$\omega$ and $Z_0, \ldots,Z_n$ are given, $Z_{n+1}$ is the sum of
$Z_n+1$ independent variables $V_{n,0},V_{n,1},\ldots,V_{n,Z_n}$
each having the geometric distribution \beqn \label{geom} P_\omega
\{V_{n,i}=j\}=\omega_{-n}(1-\omega_{-n})^j,~~~j=0,1,2,\ldots \feqn
For $s \in [0,1]$ define random variables dependent on the
environment $\omega:$ \beq n \in \zz:~~~f_n(s)&=&E_\omega
\left(s^{V_{n,0}}\right)= \fracd{1}{1+\rho_{-n}(1-s)}, \\ n \geq
0:~~~\psi_n(s)&=&E_\omega \left(s^{\sum_{i=1}^n Z_i}\right). \feq
Iterating, and taking in account that the variables $V_{n,i}$ are
conditionally independent given $Z_n$ and $\omega,$ we obtain:
\beq \psi_{n+1}(s)&=&E_\omega
\left(E_\omega\left(s^{\sum_{i=1}^{n+1} Z_i} |Z_1,Z_2,
\ldots,Z_n\right)\right)=E_\omega \left(s^{\sum_{i=1}^n Z_i}
E_\omega \left(s^{\sum_{j=0}^{Z_n} V_{n,j}} |Z_n \right)\right) =
\\&=&E_\omega \left(s^{\sum_{i=1}^n Z_i} f_n(s)^{Z_n+1}\right)=
f_n(s) E_\omega \left(s^{\sum_{i=1}^{n-1} Z_i} \bigl(s f_n(s)
\bigr)^{Z_n}\right)=
\\&=&f_n(s)\cdot f_{n-1}(s f_n(s))\cdot \ldots \cdot f_0(s f_1(s
\ldots f_{n-1}(s f_n(s)) \ldots)). \feq Thus the random variable
$\psi_{n+1}(s)$ is distributed the same as \beqn \label{phi-eq}
\varphi_{n+1}(s)=f_0(s)\cdot f_{-1}(s f_0(s))\cdot \ldots \cdot
f_{-n}(s f_{-n+1}(s \ldots f_{-1}(s f_0(s)) \ldots)). \feqn
\begin{lemma}
\label{pr1} For $n \geq 0,$ $\varphi_n(s)=1/B_n(s),$ where
$B_n(s)$ are polynomials of $s$ satisfying the following recursion
formula \beqn \label{rec1}
B_{n+1}(s)=(1+\rho_n)B_n(s)-s\rho_n B_{n-1}(s),~~~~~B_0(s)=1,~
B_1(s)=1+\rho_0(1-s). \feqn
\end{lemma}
\begin{proof}
For $n \geq 0$ let $q_n(s)=f_{-n}(s f_{-n+1}(s \ldots f_{-1}(s
f_0(s)) \ldots)).$ In particular, we have
$q_0(s)=f_0(s)=B_0(s)/B_1(s).$ By induction,
$q_n(s)=B_n(s)/B_{n+1}(s) \in (0,1]$ for all $n \in \nn.$ Indeed,
assuming that for some $n \geq 0,$ $q_{n-1}=B_{n-1}/B_n,$ we
obtain: \beq q_n=f_{-n}(s q_{n-1})=\fracd{1}{1+\rho_n \left(1-s
\fracd{B_{n-1}}{B_n}\right)}=\fracd{B_n}{ (1+\rho_n)B_n-s \rho_n
B_{n-1}}=\fracd{B_n}{B_{n+1}}. \feq It follows from \eqref{phi-eq}
that \beq \varphi_n(s)=\prod_{j=0}^{n-1}q_j(s)=
\fracd{1}{B_1(s)}\cdot \fracd{B_1(s)}{B_2(s)} \cdot \ldots \cdot
\fracd{B_{n-1}(s)}{B_n(s)}=\fracd{1}{B_n(s)}, \feq completing the
proof.
\end{proof}
Rewriting \eqref{rec1} in the form
$B_{n+1}-B_n=\rho_n(B_n-B_{n-1})+\rho_n(1-s)B_{n-1}$ and
iterating, we obtain another, useful in the sequel, form of the
recursion: \beqn \label{core} B_{n+1}(s)=B_n(s)+(1-s) \cdot
\suml_{i=0}^n B_{i-1}(s) \prodl_{j=i}^n \rho_j,~~~~~B_{-1}(s)=
B_0(s)=1. \feqn In particular, $B_n(s)$ is a strictly increasing
sequence for a fixed $s \in [0,1).$
\begin{proposition}
\label{scale-zed} Suppose that Assumption \ref{v239} holds. Then,
\beqn \label{scz} n^{-1/\alpha} \suml_{i=1}^n Z_i
\overset{\pp}{\Ra}
\begin{cases}
0&\mbox{\em if}~\alpha<\kappa\\
\infty&\mbox{\em if}~\alpha>\kappa,
\end{cases}
\feqn where $\kappa$ is defined in \eqref{kappa}.
\end{proposition}
\begin{proof}
For $m \geq 0$ let \beqn \label{defini}
H_m(\omega)&=&\rho_m+\rho_m\rho_{m-1}+...+\rho_m\rho_{m-1}\ldots
\rho_0. \nonumber \feqn First, suppose that $0<\alpha<\kappa.$
Choose any $\beta \in (\alpha,\kappa).$ By Chebyshev's inequality,
we obtain that for any $x>0,$ \beqn \label{cheb} \pp \left(
\fracd{1}{n^{1/\alpha}} \sum_{i=1}^n Z_i >x \right) \leq
\fracd{1}{x^\beta n^{\beta/\alpha}} \ee \left( \left( \sum_{i=1}^n
Z_i \right)^{\beta} \right) .\feqn Using the inequality $\left(
\sum_{i=1}^n Z_i \right)^\beta \leq \sum_{i=1}^n Z_i ^\beta$ and
Jensen's inequality $E_\omega\bigl(Z_i^\beta \bigr) \leq
\bigl(E_\omega(Z_i) \bigr)^\beta,$ \beq \ee \left( \left(
\sum_{i=1}^n Z_i \right)^\beta \right)& \leq& \ee \left(
\sum_{i=1}^n Z_i^\beta \right)=\sum_{i=1}^n E_P \left(E_\omega
\left( Z_i^\beta \right) \right) \leq \sum_{i=1}^n E_P \left(
\left(E_\omega \left( Z_i\right)\right)^\beta \right) =
\\&=&\sum_{i=1}^n E_P \left( H_{i-1}^\beta \right) \leq n E_P \left(
R^\beta \right) ,\feq where $R$ is defined in \eqref{sf}. Since by
\eqref{sign}, $ E_P\left( R^\beta \right)<\infty$ for all
$\beta<\kappa,$ the claim for $\alpha \in (0,\kappa)$ follows from
\eqref{cheb}.
\\
\centerline{} \\
Now suppose that $\alpha>\kappa.$ Fix a real number $\lambda>0$
and let $s_n=e^{-\lambda/n^{1/\alpha}}.$ In order to prove that
\eqref{scz} holds for $\alpha>\kappa,$ it is sufficient to show
that $P\left(\lim_{n \to \infty} B_n(s_{n+1}) =\infty \right)=1$
for any $\lambda \geq 0.$ We next prove that for any $M \in \rr,$
$P\bigl(B_n(s_{n+1})<M~\mbox{i.o.}\bigr)=0,$ implying the latter
assertion.

If $B_n(s_{n+1})< M,$ then $1 \leq B_0(s_{n+1}) \leq
B_1(s_{n+1})\leq \ldots \leq B_n(s_{n+1})< M.$ It follows from
(\ref{core}) that for $m \leq n ,$ \beq
1&=&\fracd{B_{m-1}(s_{n+1})}{B_m(s_{n+1})}+(1-s_{n+1})
\suml_{i=0}^{m-1} \fracd{B_{i-1}
(s_{n+1})}{B_m(s_{n+1})}\prodl_{j=i}^{m-1} \rho_j \geq
\fracd{B_{m-1}(s_{n+1})}{B_m(s_{n+1})}+\fracd{1}{M}(1-s_{n+1})H_{m-1}.
\feq Thus, $B_n(s_{n+1})< M$ implies that $H_m \leq M \lambda^{-1}
(n+1)^{1/\alpha},$ for $m=0,1,\ldots,n-1.$ We conclude that if
$\kappa< \gamma < \alpha $ and $n$ large enough, then
$B_n(s_{n+1})< M $ implies \beqn \label{io} H_m \leq
(n+1)^{1/\gamma},~~~~~m=0,1,\ldots,n-1. \feqn Let $K_n=[1/\delta
\cdot \log n],$ where $[\hspace{2pt} \cdot \hspace{2pt}]$ denotes
the integer part of a number and $\delta>0$ is a parameter whose
precise value will be determined later. Further, let
$\veps=\fracd{1}{3}\left(1-\kappa/\gamma\right),$
$L_n=K_n+n^\veps$ and \beq G_{j,n}=\prodl_{i=(j-1)K_n}^{j K_n-1}
\rho_i ,~~j=1,2,\ldots \left[n/K_n\right]. \feq From the LDP for
$\ol{R}_n=\fracd{1}{n} \suml_{i=0}^{n-1} \log \rho_i$ we obtain
that there is $N_\veps>0$ such that $n>N_\veps$ implies \beqn
\label{est} P\bigl( G_{1,n} >(n+1)^{1/\gamma} \bigr)> n^{-1/\delta
\cdot J \left(\delta /\gamma\right)-\veps}. \feqn It follows from
(\ref{io}) and (\ref{est}) that for $n$ large enough, \beq
P\bigl(B_n(s_{n+1})< M\bigr)&\leq&P\Bigl( G_{i L_n,n}
<(n+1)^{1/\gamma},~~i=1,2,\ldots [n/L_n]-1 \Bigr) \leq
\\&\leq&
\left(1-n^{-1/\delta \cdot J \left(\delta /\gamma
\right)-\veps}\right)^{[n/L_n]}+ [n/L_n] \cdot
\eta\left(n^\veps\right), \feq where $\eta(l)$ are the strong
mixing coefficients. That is, for $n$ large enough, \beqn
\label{ozenka} P(B_n(s_{n+1})< M) \leq \exp
\left\{-\fracd{1}{2}n^{1-1/\delta \cdot J
    \left(\delta /\gamma \right)- 2\veps} \right \}+n^{1-\veps} \cdot
\eta\left(n^\veps \right).  \feqn Now, we will choose the
parameter $\delta$ in such a way that $1-1/\delta \cdot J(\delta
/\gamma)-2\veps =\veps.$ By \eqref{kappa}, there exists $x^*>0$
such that $J(x^*)=\kappa x^*.$ Hence, letting $\delta=\gamma x^*$
we obtain, \beq 1-1/\delta \cdot J(\delta /\gamma)-2\veps
=1-\kappa/\gamma -2\veps=\veps.  \feq It follows from
(\ref{ozenka}) that for $n$ large enough, \beq
P\bigl(B_n(s_{n+1})< M \bigr) \leq e^{-n^\veps/2}+n^{1-\veps}
\cdot \eta\left(n^\veps \right).  \feq Thus, $P\bigl(B_n(s_{n+1})<
M~\mbox{i.o.}\bigr)=0$ by the Borel-Cantelli lemma.
\end{proof}
Proposition \ref{scale-ti} follows from \eqref{ti-zed},
\eqref{ti-minus}, and Proposition \ref{scale-zed}.
\section*{Acknowledgments}
I would like to thank my thesis advisors Eddy Mayer-Wolf and Ofer
Zeitouni for posing the problem studied in this paper and for many
helpful suggestions.

\providecommand{\bysame}{\leavevmode\hbox
to3em{\hrulefill}\thinspace}
\providecommand{\MR}{\relax\ifhmode\unskip\space\fi MR }
\providecommand{\MRhref}[2]{%
  \href{http://www.ams.org/mathscinet-getitem?mr=#1}{#2}
} \providecommand{\href}[2]{#2}

 \end{document}